\input amstex
\input amsppt.sty

\magnification1200
\vsize19cm


\def\LL{\leavevmode\setbox0=\hbox{L}\hbox to\wd0{\hss\char'40L}}
\def\al{\alpha}

\def\ep{\varepsilon}

\def\la{\lambda}
\def\rh{\rho}

\def\ta{\sigma}

\def\om{\omega}


\def\R{{\Bbb R}}
\def\S{{\Bbb S}}

\def\C{{\Bbb C}}

\def\A{{\Cal A}}

\def\N{{\Bbb N}}
\def\Z{{\Bbb Z}}

\def\today{\ifcase\month\or
 January\or February\or March\or April\or May\or June\or
 July\or August\or September\or October\or November\or December\fi
 \space\number\day, \number\year}

\def\({\left(}
\def\){\right)}
\def\[{\left[}
\def\]{\right]}

\def\suml{\sum\limits}

\def\3{\ss}

\topmatter
\title
Some conjectures for\\ Macdonald polynomials of type B, 
C, D
\endtitle 
\rightheadtext{Macdonald polynomials}
\author Michel Lassalle 
\endauthor 
\affil 
Centre National de la Recherche Scientifique\\
Institut Gaspard Monge, Universit\'e de Marne-la-Vall\'ee\\
77454 Marne-la-Vall\'ee Cedex, France
\endaffil
\dedicatory
Dedicated to Alain Lascoux on the occasion of his 60th birthday
\enddedicatory
\address 
Institut Gaspard Monge, Universit\'e de Marne-la-Vall\'ee, 
77454 Marne-la-Vall\'ee Cedex, France
\endaddress
\email
lassalle\@univ-mlv.fr
\endemail
\urladdr
http://igm.univ-mlv.fr/${}^{\sim}$lassalle
\endurladdr

\abstract  We present conjectures giving formulas 
for the Macdonald polynomials of type $B$, $C$, $D$ 
which are indexed by a multiple of the first 
fundamental weight. The transition matrices 
between two different types are explicitly given. 
\endabstract
\endtopmatter

\document

\parskip.1cm plus 3pt minus -3pt

\head \bf Introduction\endhead

Among symmetric functions, the special importance of Schur 
functions comes from their intimate connection with 
representation theory. Actually the irreducible polynomial 
representations of $GL_{n}(\C)$ are indexed by partitions 
$\la=(\la_{1},\ldots,\la_{n})$ of length $\le n$, and their 
characters are the Schur functions $s_{\la}$. 

In the eighties, I. G. Macdonald introduced a new family of symmetric 
functions $P_{\la}(q,t)$. These orthogonal polynomials depend rationally on two 
parameters $q,t$ and generalize Schur functions, which are 
obtained for $t=q$ \cite{9,10}. 

When the indexing partition is reduced to a row $(k)$ (i.e. has 
length one), the Macdonald polynomial $g_k(q,t)$ of $n$ 
variables $x=(x_1,\ldots,x_n)$ are given by their generating function
$$\prod_{i=1}^n \frac{(tux_i;q)_\infty} {(ux_i;q)_\infty}\,
=\sum_{r \ge 0}u^r g_r(x;q,t),$$
with the standard notation $(a;q)_\infty=\prod _{i=0}^{\infty}(1-aq^i)$.
Of course for $t=q$ the complete functions $s_{(r)}=h_r$ are recovered.

A few years later, generalizing his previous work, 
Macdonald introduced another class of orthogonal polynomials, 
which are Laurent polynomials in several variables, and generalize 
the Weyl characters of compact simple Lie groups \cite{11,12}.
In the most simple situation of this new framework, a family 
$P_{\la}^{(R)}(q,t)$ of polynomials, depending rationally 
on two parameters $q,t$, is attached to each root system $R$. 

These orthogonal polynomials are elements of the group algebra of 
the weight lattice of $R$, invariant under the action of the Weyl 
group. They are indexed by the dominant weights of $R$.  

When $R$ is of type $A$, the orthogonal polynomials $P_{\la}^{(R)}(q,t)$ 
correspond to the symmetric functions $P_{\la}(q,t)$ previously studied 
in \cite{9,10}. For $t=q$, they  correspond to the Weyl characters 
$\chi_{\la}^{(R)}$ of compact simple Lie groups.

This paper is only devoted to the Macdonald polynomials 
which are indexed by a multiple of the first fundamental 
weight $\om_1$. Since H. Weyl \cite{15}, it is well known that 
$\chi_{r\om_1}^{(R)}$ is given by
\smallskip\noindent
(i) $h_{r}(X)+h_{r-1}(X)$, when $R=B_n$,

\smallskip\noindent
(ii) $h_{r}(X)$, when $R=C_n$,

\smallskip\noindent
(iii) $h_{r}(X)-h_{r-2}(X)$, when $R=D_n$,

\noindent
with $X=(x_1,\ldots,x_n,1/x_1,\ldots,1/x_n)$.
 
However, as far as the author is aware, no such result is 
known when $t \neq q$, and no explicit expansion is 
available for the Macdonald polynomials 
$P_{r\om_1}^{(R)} (q,t)$. The purpose of this paper is to present 
some conjectures generalizing the previous formulas.

Actually this problem can be considered in a more general 
setting, allowing two distinct parameters $t,T$, 
each of which is attached to a length of roots. We give an 
explicit formula for $P_{r\om_1}^{(R)} (q,t,T)$ when $R$ is of 
type $B,C,D$, together with an explicit formula for the transition 
matrices between different types. The entries of these transition 
matrices appear to be fully factorized and reveal deep connections with basic hypergeometric series. 

We provide some support for these conjectures by showing that they are verified upon principal specialization. On the other hand, computer calculations show a very strong empirical evidence in their favor. 

\head \bf 1. Macdonald polynomials \endhead

In this section we introduce our notations, and recall some 
general facts about Macdonald polynomials. For more details the reader is referred to \cite{11,12,13}. 

The most general class of Macdonald polynomials is associated with a 
pair of root systems $(R,S)$, spanning the same vector space and 
having the same Weyl group, with $R$ reduced. Here we shall only 
consider the case of a pair $(R,R)$, with $R$ of type $B,C,D$.

Let $V$ be a finite-dimensional real vector space endowed with a
positive definite symmetric bilinear form $\langle  u, v 
\rangle$. For all $v\in V$, we write $\vert v\vert = \langle v, v \rangle ^{1/2}$ , and $v^\lor =2v/\vert v\vert^2$.

Let $R\subset V$ be a reduced irreducible root system, $W$ the Weyl group of $R$, $R^+$ the set of positive roots, $\{\alpha_1,\dots,\alpha_n\}$ 
the basis of simple roots, and $R^\lor  =\{ \alpha^\lor \mid \alpha\in R\}$ the dual root system.

Let $Q=\sum_{i=1}^n\Z\ \alpha_i$ and 
$Q^+=\sum_{i=1}^n\N\ \alpha_i$ be the root lattice 
of $R$ and its positive octant. Let $P=\{\la \in V \ |\
\langle \la , \alpha^\lor \rangle \in \Z \ \forall\al\in R \}$ and 
$P^{+}=\{\la \in V \ |\ \langle \la , \al^\lor \rangle \in 
\N \ \forall\al\in R^+ \}$ 
be the weight lattice of $R$ and the cone of dominant weights.

A basis of $Q$ is formed by the simple roots $\alpha_i$. A basis of $P$ is formed by the fundamental weights $\omega_i$ 
defined by $\langle\omega_{i}, \alpha_{j}^\lor \rangle=\delta _{ij}$. 
A partial order is defined on $P$ by
$\lambda\geq\mu$ if and only $\lambda-\mu\in Q^+$.

Let $\A$ denote the group algebra over $\R$ of the free Abelian group
$P$. For each $\lambda\in P$ let $e^\lambda$ denote the corresponding
element of $\A$, subject to the multiplication rule 
$e^\lambda e^\mu = e^{\lambda +\mu}$. The set $\{e^\lambda, \lambda\in P\}$ forms an $\R$-basis of $A$.

The Weyl group $W$ acts on $P$ and on $\A$.
Let $\A^{W}$ denote the subspace of $W$-invariants in $A$. 
Such elements are called ``symmetric polynomials''.
There are two important examples of a basis of $\A^{W}$. The 
first one is given by the orbit-sums
$$m_\la = \suml_{\mu\in W\la}e^{\mu}, \quad \quad \quad \la \in
P^{+}.$$
Another basis is provided by the Weyl characters defined as 
follows. Let
$$\delta =\prod_{\alpha\in R^+}
(e^{\alpha/2}-e^{-\alpha/2})=e^{-\rh} \prod_{\alpha\in
R^+} (e^{\al}-1),$$
with $\rh=\frac{1}{2} \sum_{\alpha \in R^{+}} \alpha \in P$.
Then $w\delta = \varepsilon(w)\delta$ for any $w\in W$, where
$\varepsilon(w)=\det(w)=\pm 1$. For all $\lambda\in P$,
$$\chi_\lambda = \delta^{-1}\suml_{w\in W}\varepsilon(w)e^{w(\lambda +
\rh)}$$
is in $\A^{W}$, and the set $\{\chi_\lambda,\ \lambda\in P^{+}\}$ 
forms an $\R$-basis of $\A^{W}$.

Let $0< q <1$. For any indeterminate $x$ and 
for all $k\in\N$, define
$$(x;q)_\infty=\prod _{i=0} ^{\infty}(1-xq^i),
\quad \quad (x;q)_k=\prod _{i=0} ^{k-1}(1-x q^i).$$
For each $\alpha\in R$ let $t_\alpha=q^{k_\al}$ be a positive real 
number such that $t_\alpha=t_\beta$ if $\vert\alpha\vert=\vert\beta\vert$. 
Then we have at most two different values for the $t_\alpha$'s. 
Define
$$\rh_k=\frac{1}{2} \suml_{\alpha \in R^{+}} k_\al \al,
\quad \quad \rh_k^\ast=\frac{1}{2} \suml_{\alpha \in R^{+}} k_\al \al^\lor.$$

If $f=\sum_{\la \in P} \ a_{\la} e^{\la} \in \A$, let
$\overline{f}=\sum_{\la \in P} \ a_{\la}e^{-\la }$ and
$[\,f\,]_{1}$ its constant term $a_{0}$.
The inner product defined on $\A$ by
$$\langle f, g\rangle _{q,t}=\frac{1}{|W|}[f\bar g\Delta_{q,t}]_{1},$$ 
with $|W|$ the order of $W$, and
$$\Delta_{q,t} =\prod_{{\alpha\in R}}\frac 
{(e^\alpha;q)_\infty}{(t_\alpha e ^\alpha;q)_\infty}$$
is non degenerate and $W$-invariant.

There exists a unique basis $\{P_\la, \ \la \in P^{+}\}$ of $\A^W$, 
called Macdonald polynomials, such that
\smallskip\noindent
(i) $P_\la = m_\la +
\sum_{\mu\in P^+,\ \mu < \la} \ a_{\la\mu}(q,t)\ m_\mu$\newline
where the coefficients $a_{\la\mu}(q,t)$ are rational functions of 
$q$ and the $t_\alpha$'s,

\smallskip\noindent
(ii) $\langle P_\la, P_\mu\rangle  _{q,t}= 0$ if $\la\ne\mu$.

It is clear that the $P_\la$, if they exist, are unique. Their 
existence is proved as eigenvectors of an operator $E:\A^W\rightarrow \A^W$, 
selfadjoint with respect to $\langle \ ,\ \rangle _{q,t}$, 
and having its eigenvalues all distinct. When $R$ is of 
type $B,C,D$ this operator may be constructed as follows 
\cite{11,12}.

Let $\pi$ be a minuscule weight of $R^\lor$, 
i.e. a vector $\pi\in  V$ such
that $\langle \pi,\alpha \rangle$ takes only values $0$ and
$1$ for $\alpha \in R^+$. Such a vector exists when $R$ is of 
type $B,C,D$, and is necessarily a fundamental weight of $R^\lor $. Let
$$\Phi_\pi =
\prod_{\alpha\in R^+}\frac {1-t_\alpha^{\langle\pi,\alpha\rangle} e^\alpha} 
{1-e^\alpha},$$ 
and $T_\pi$ the translation operator defined on $\A$ by
$T_{\pi }(e^{\la })=q^{\langle\la ,\pi \rangle}e^{\la }$ for any $\la \in P$.
Let $E_\pi$ the operator defined by
$$E_\pi f = \sum_{w\in W} w(\Phi_\pi\cdot T_\pi f).$$
Macdonald polynomials $P_\la$ are then introduced as eigenvectors 
of $E_\pi$ :
$$E_\pi\ P_\la=c_{\la}\ P_\la
\quad \quad \text{with} \quad \quad
c_{\la}=q^{\langle\pi,\rh _k\rangle}\sum_{w\in W}q
^{\langle w\pi,\lambda + \rh_k\rangle}.$$

We may regard any $f=\sum_{\la \in P}  f_\lambda e ^\lambda \in \A$,
having only finitely many nonzero coefficients, as a function on $V$ by putting 
for any $x\in V$,
$$f(x)=\sum_{\la \in P} f_\lambda q ^{\langle\lambda,x\rangle}.$$
Then Macdonald polynomials satisfy the following 
Specialization Formula \cite{2}
$$P_\la(\rh_k^\ast)= q^{-\langle\la,\rh_k^\ast\rangle} \,
\prod_{\al \in R^+} 
\frac{(q^{\langle\rh_k,\al^\lor\rangle}t_{\al};q)_{\langle\la,\al^\lor\rangle}}
{(q^{\langle\rh_k,\al^\lor\rangle};q)_{\langle\la,\al^\lor\rangle}}.$$

In the sequel we shall 
consider $R$ to be of type $B_n, C_n$ or $D_n$. We shall
identify $V$ with $\R^n$ and write $\ep_1,\ldots,\ep_n$ for 
its standard basis. 
Defining $x_i=e^{\ep_i}$ for $i=1\ldots n$, we shall regard 
$P_\la$ as a Laurent polynomial of $n$ variables $x_1,\ldots,x_n$.

We shall only consider Macdonald polynomials $P_\la$ 
associated with a weight $\la=r\om_1$, multiple of 
the first fundamental weight $\om_1=\ep_1$.

\head \bf 2. Basic hypergeometric series \endhead

We shall need three results about basic hypergeometric series. 
The author is deeply indebted to Professor Mizan Rahman for 
communicating their proofs to him. Since these results have intrinsic interest, their proofs are given below in Section 10.

We adopt the notation of \cite{3} and write
$${}_{r+1}\phi_r \left[\matrix
a_1,a_2,\dots,a_{r+1}\\
b_1,b_2,\dots,b_r \endmatrix;q,z \right]=
\sum_{i\ge 0}\frac{(a_1;q)_i\ldots(a_{r+1};q)_i}
{(b_1;q)_i\ldots(b_r;q)_i}\frac{z^i}{(q;q)_i}.$$

\proclaim{Theorem 1}We have the following transformation 
between ${}_2{\phi}_1$ series
$$\multline
\frac{(u;q)_r}{(q;q)_r}\ {}_2{\phi}_1 
\left[\matrix q^{-r}, \ ux\\q^{1-r}/ux \endmatrix;q,qv/u^2x \right]=
\sum_{i=0}^{[r/2]} 
\frac{(u;q)_{r-2i}}{(q;q)_{r-2i}}\ 
{}_2{\phi}_1 
\left[\matrix q^{2i-r}, \ uy\\q^{1+2i-r}/uy \endmatrix;q,qv/u^2y \right]\\
\times y^iv^i\ \frac{(x/y;q)_i}{(q;q)_i} \
\frac{(uq^{r-2i};q)_{2i}}{(uxq^{r-i};q)_{i}\ (uyq^{r-2i+1};q)_{i}}.
\endmultline$$
\endproclaim

An infinite-dimensional matrix $(f_{ij})_{i,j\in \Z}$ is said to 
be lower-triangular if $f_{ij}=0$ unless $i\ge j$. 
Two infinite-dimensional lower-triangular 
matrixes $(f_{ij})_{i,j\in \Z}$ and $(g_{kl})_{k,l\in \Z}$ 
are said to be mutually inverse if 
$\sum_{i\ge j \ge k} f_{ij}g_{jk}=\delta_{ik}$.

\proclaim{Corollary} Defining
$$\Cal M^{(u,v;x,y)}_{r,r-2i}=y^iv^i\ \frac{(x/y;q)_i}{(q;q)_i} \
\frac{(uq^{r-2i};q)_{2i}}{(uxq^{r-i};q)_{i}\ (uyq^{r-2i+1};q)_{i}},$$
the infinite matrices $\Cal M(u,v;x,y)$ and $\Cal M(u,v;y,x)$ 
are mutually inverse.
\endproclaim

Michael Schlosser remarked that this corollary is equivalent 
with Bressoud's matrix inverse~\cite{1}, which states that, defining
$$\Cal A^{(u,v)}_{ij}=(u/v)^j\ \frac{(u/v;q)_{i-j}}{(q;q)_{i-j}} \
\frac{(u;q)_{i+j}}{(vq;q)_{i+j}}\ 
\frac{1-vq^{2j}}{1-v},$$
the matrices $\Cal A(u,v)$ and $\Cal A(v,u)$ 
are mutually inverse. Indeed let $d$ be either $0$ or $1$. Replacing $r$ by 
$2r+d$, and $i$ by $r-i$, some factors cancel or can be
pulled out of the previous sum, yielding the above form of~\cite{1}.

Bressoud's matrix inverse was originally derived
from the terminating very-well-poised $_6\phi_5$
summation~\cite{3, (II.21)}. We refer to~\cite{4, 8} for some generalizations 
of~\cite{1}, as well as more details and references about inversion 
of infinite matrices. 

\proclaim{Theorem 2} Defining
$$\multline
\frac{(u;q)_{r}}{(q;q)_{r}}\ H_r=
\sum_{i=0}^r a^{r-i}\ \frac{(u;q)_{r-i}}{(q;q)_{r-i}}\
{}_2{\phi}_1 
\left[\matrix q^{i-r},\ v\\q^{1+i-r}/v 
\endmatrix;q,q/v^2a^2 \right]\\\times 
\frac{(x;q)_{i}}{(q;q)_{i}}\
\frac{(uq^{r-i};q)_{i}}{(vq^{r-i+1};q)_{i}}\
\frac{(v^2q^{2r-i+1};q)_{i}}{(xv^2q^{2r-i};q)_{i}},
\endmultline$$
we have
$$H_r=(-1/v)^r\frac{(xv^2;q)_r}{(xv^2;q)_{2r}}\
 \frac{(x^2v^2;q^2)_r(v^2q;q^2)_r}{(xv;q)_{r}}\
{}_4{\phi}_3 \left[\matrix q^{-r},xv^2q^r,-av,-1/a\\
vq^{\frac{1}{2}},-vq^{\frac{1}{2}},-xv\endmatrix;q,q \right].$$
\endproclaim

A converse property may be stated as follows.
\proclaim{Theorem 3} 
We have
$$\multline
a^r \frac{(u;q)_{r}}{(q;q)_{r}}\
{}_2{\phi}_1 
\left[\matrix q^{-r},\ v\\q^{1-r}/v 
\endmatrix;q,q/v^2a^2 \right]
= \sum_{i=0}^r \frac{(u;q)_{r-i}}{(q;q)_{r-i}}\ H_{r-i}
\\\times x^i\ \frac{(1/x;q)_{i}}{(q;q)_{i}}\
\frac{(uq^{r-i};q)_{i}}{(vq^{r-i+1};q)_{i}}\
\frac{(v^2q^{2r-2i+1};q)_i}
{(xv^2q^{2r-2i+1};q)_{i}}\
\frac{1-v^2q^{2r}}{1-v^2q^{2r-i}}.
\endmultline$$
\endproclaim

Observe that, as a consequence of Theorems 2 and 3, we recover the special case $y=1$ of the following matrix inverse.
\proclaim{Lemma 1} Defining 
$$\Cal N^{(u,v;x,y)}_{ij}=y^{i-j} \
\frac{(x/y;q)_{i-j}}{(q;q)_{i-j}}\
\frac{(uq^j;q)_{i-j}}{(vq^{j+1};q)_{i-j}}\
\frac{(v^2q^{2j+1};q)_{2i-2j}}
{(xv^2q^{i+j};q)_{i-j}\ (yv^2q^{2j+1};q)_{i-j}},$$
the infinite matrices $\Cal N(u,v;x,y)$ and $\Cal N(u,v;y,x)$ 
are mutually inverse.
\endproclaim
\demo{Proof}
Defining
$$\Cal B^{(x,y)}_{ij}=y^{i-j}\ 
\frac{(x/y;q)_{i-j}}{(q;q)_{i-j}} \
\frac{1}{(xq^{i+j};q)_{i-j}\ (yq^{2j+1};q)_{i-j}},$$
the infinite matrices 
$\Cal B(x,y)$ and $\Cal B(y,x)$ are mutually inverse, as a consequence of a result of Krattenthaler \cite{4}. 
If two infinite matrices $(f_{ij})$ and $(g_{kl})$ 
are mutually inverse, for any sequence $(d_k)$, the 
matrices $(f_{ij}\ d_{i}/d_{j})$ and $(g_{kl}\ d_{k}/d_{l})$ 
are obviously mutually inverse. 
Since
$$(v^2q^{2j+1};q)_{2i-2j}\ \frac{(uq^j;q)_{i-j}}{(vq^{j+1};q)_{i-j}}=
\frac{(v^2q;q)_{2i}}{(v^2q;q)_{2j}}\
\frac{(u;q)_i}{(u;q)_j}\ \frac{(vq;q)_j}{(vq;q)_i},$$
we apply this property to $\Cal B(xv^2,yv^2)$ and 
$\Cal B(yv^2,xv^2)$, with
$$d_k=v^{2k}\ (v^2q;q)_{2k}\ \frac{(u;q)_k}{(vq;q)_k}.
\quad \quad \qed$$
\enddemo

\head \bf 3. Type $C$ \endhead

For the type $C$ root system, the set of positive roots is the union 
of $R_{1}=\{\ep_{i}\pm \ep_{j}, 1\le i < j\le n\}$ and 
$R_{2}=\{2\ep_{i}, 1\le i\le n\}$. Elements of each 
set have the same length, and we write $t_{\al}=t$ for $\al \in R_{1}$ 
and $t_{\al}=T$ for $\al \in R_{2}$. 

The Weyl group $W$ is the semi-direct product of the 
permutation group $S_{n}$ by $(\Z/2\Z)^n$. 
It acts on $V$ by signed permutation of components. 
The fundamental weights are given by 
$\om_{i}=\sum_{j=1}^i\ep_j,\ 1\le i\le n$. The dominant 
weights $\la\in P^+$ can be identified with vectors 
$\la=\sum_{i=1}^n \la_i \ep_i$, such that
$(\la_1,\la_2,\ldots,\la_n)$ is a partition. 
There is only one minuscule 
weight $\om_1$. The partial order $\la\geq\mu$ is given by
$$\sum_{i=1}^j (\la_i-\mu_i) \in \N, \quad \text{for} 
\quad j=1,\ldots, n-1, \quad \quad
\sum_{i=1}^n (\la_i-\mu_i) \in 2\N.$$
With $t=q^{k}$ and $T=q^{K}$, we have 
$$\rh_k=\sum_{i=1}^n \left((n-i)k+K \right) \ep_i,\quad \quad
  \rh_k^\ast= \sum_{i=1}^n \left((n-i)k+K/2 \right) \ep_i.$$

If we write $x_i=e^{\ep_i}$ for $i=1\ldots n$, and define
$$P_{r\om_1}^{(C)}(x;q,t,T)=\frac{{(q;q)_r}}{{(t;q)_r}} 
\ g_r^{(C)}(x;q,t,T),$$
the Specialization Formula reads
$$g_r^{(C)}(t^{n-1} T^{\frac{1}{2}},t^{n-2} 
T^{\frac{1}{2}},\ldots, T^{\frac{1}{2}} ;\ q,t,T)=
\frac{t^{r(1-n)}}{T^{r/2}} \ \frac{(t^n;q)_r}{(q;q)_r}\
\frac{(T^2t^{2n-2};q)_r}{(Tt^{n-1};q)_r}.$$

Let us introduce the 
auxiliary quantities $G_r(x;q,t)$ defined by their generating function
$$\prod_{i=1}^n
\frac{(tux_i;q)_\infty} {(ux_i;q)_\infty}\
\frac{(tu/x_i;q)_\infty} {(u/x_i;q)_\infty}
=\sum_{r \ge 0}u^r G_r(x;q,t).$$
In $\la$ - ring notation (see Section 8), they can be written as 
$$G_r(x;q,t)=h_r\left[\frac{1-t}{1-q}X^{\dag}\right],$$
with $X^{\dag}=\sum_{i=1}^n (x_i+1/x_i)$. Their 
specialization may be given as follows.

\proclaim{Lemma 2} For any positive integer $r$ we have
$$G_r(t^{n-1}a,t^{n-2}a,\ldots,a ;\ q,t)
=a^{r}\ \frac{(t^n;q)_r}{(q;q)_r}\ {}_2{\phi}_1 
\left[\matrix q^{-r},\ t^n\\t^{-n}q^{1-r} \endmatrix;q,qt^{1-2n}/a^2 \right].$$
\endproclaim
\demo{Proof}
Taking into account
$$\prod_{i=1}^n
\frac{(t^{n-i+1}ua;q)_\infty} {(t^{n-i}ua;q)_\infty}\
\frac{(t^{i-n+1}u/a;q)_\infty} {(t^{i-n}u/a;q)_\infty}
=\frac{(t^nua;q)_\infty}{(ua;q)_\infty}\
\frac{(tu/a;q)_\infty}{(t^{1-n}u/a;q)_\infty},$$
and applying the classical $q$--binomial formula
\cite{3, (II.3)}
$$\frac{(tu;q)_\infty} {(u;q)_\infty}
=\sum_{i \ge 0} \frac{(t;q)_i}{(q;q)_i} \ u^i,$$
we have
$$G_r(t^{n-1}a,t^{n-2}a,\ldots,a ;\ q,t)=
\sum_{i=0}^r a^{(r-2i)} t^{i(1-n)}\ \frac{(t^n;q)_i}{(q;q)_i} \
\frac{(t^n;q)_{r-i}}{(q;q)_{r-i}}.\qed$$
\enddemo

\proclaim{Conjecture 1} For any positive integer $r$ we have
$$g_r^{(C)}(q,t,T)=\sum_{i=0}^{[r/2]}
G_{r-2i}(q,t)\ t^i \ \frac{(T/t;q)_i}{(q;q)_i} \
\frac{(t^nq^{r-i};q)_{i}}{(Tt^{n-1}q^{r-i};q)_{i}}\
\frac{1-t^nq^{r-2i}}{1-t^nq^{r-i}}.$$
Conversely
$$G_r(q,t)=\sum_{i=0}^{[r/2]}
g_{r-2i}^{(C)}(q,t,T)\  T^i \ \frac{(t/T;q)_i}{(q;q)_i} \
\frac{(t^nq^{r-2i};q)_{i}}{(Tt^{n-1}q^{r-2i+1};q)_{i}}.$$
In other words, the transition matrix from $g^{(C)}(q,t,T)$ to 
$G(q,t)$ is $\Cal M (t^n,t;T/t,1)$, and 
its inverse is $\Cal M (t^n,t;1,T/t)$. 
\endproclaim

Using $\la$ - ring techniques, we have proved this conjecture for $T=t$. 
\proclaim{Theorem 4}
For any positive integer $r$ we have $g_r^{(C)}(q,t,t)=G_r(q,t)$.
\endproclaim

The proof will be given below in Section $9$. The Specialization Formula gives some support to Conjecture $1$. 

\proclaim{Lemma 3} Conjecture 1 yields the specialization
$$g_r^{(C)}(t^{n-1} a,t^{n-2} a,\ldots, a ;\ q,t,T) = a^r\ 
\frac{(t^n;q)_r}{(q;q)_r}\ {}_2{\phi}_1 \left[\matrix  q^{-r},\ Tt^{n-1} \\
\ t^{1-n}q^{1-r}/T \endmatrix;q,qt^{2-2n}/Ta^2 \right].$$
\endproclaim

\demo{Proof} A straightforward application of Theorem $1$ 
and Lemma 2 with $u=t^n$, $v=t/a^2$, $x=T/t$, and $y=1$.\qed
\enddemo

\proclaim{Corollary} Conjecture 1 is true for 
$x_i=t^{n-i} T^{\frac{1}{2}}$ $(1\le i\le n)$.
\endproclaim
\demo{Proof} Keeping the same notations, we now have $v=1/x$. 
The previous result is a $q$--Chu--Vandermonde sum \cite{3, (II.7)} given by
$$\align
T^{r/2} \frac{(u;q)_r}{(q;q)_r}\ {}_2{\phi}_1 \left[\matrix  
q^{-r},\ uv \\\ q^{1-r}/uv \endmatrix;q,q/u^2v^2 \right]&
=T^{r/2} \frac{(u;q)_r}{(q;q)_r}\
\frac{(q^{1-r}/u^2v^2;q)_r}{(q^{1-r}/uv;q)_r}\\ 
&=T^{r/2}(uv)^{-r}\ \frac{(u;q)_r}{(uv;q)_r}\
\frac{(u^2v^2;q)_r}{(q;q)_r}.
\endalign$$
We recover the Specialization Formula. \qed
\enddemo

\head \bf 4. $D$ versus $C$ \endhead

The root system $D_n$ is self-dual: $R=R^\lor$. The set of positive roots 
is $R^+=\{\ep_{i}\pm \ep_{j}, 1\le i < j\le n\}$. Roots have all
the same length and we write  $t_{\al}=t$ for $\al \in R^+$.

The Weyl group $W$ is the semi-direct product of the 
permutation group $S_{n}$ by $(\Z/2\Z)^{n-1}$. 
It acts on $V$ by signed permutation of components, subject 
to the condition that the number of minus signs is even. 
The fundamental weights are given by 
$\om_{i}=\sum_{j=1}^i\ep_j$ for $1\le i\le n-2$. The
``spin weights'' $\om_{n-1}$ and  $\om_{n}$ are  defined by
$\om_{n}=\frac{1}{2}(\ep_1+\ldots+\ep_n)$ and 
$\om_{n-1}=\om_{n}-\ep_n$. There are three minuscule 
weights $\om_1$, $\om_{n-1}$ and $\om_{n}$.

The dominant weights $\la\in P^+$ can be identified with vectors 
$\la=\sum_{i=1}^n \la_i \ep_i$, whose components are all integers or all 
half-integers, and subject to the condition
$\la_1 \ge \la_2\ge \ldots\ge \la_{n-1}\ge|\la_n|$. 
The partial order $\la\geq\mu$ is given by
$$\sum_{i=1}^j (\la_i-\mu_i) \in \N, \quad \text{for} 
\quad j=1,\ldots, n-2, \quad \quad
\sum_{i=1}^{n-1} (\la_i-\mu_i) \ \pm (\la_n-\mu_n) \in 2\N.$$
Writing $t=q^{k}$ we have 
$\rh_k=\rh_k^\ast=k \ \suml_{i=1}^n (n-i) \ep_i.$

If we define
$$P_{r\om_1}^{(D)}(x;q,t)=\frac{{(q;q)_r}}{{(t;q)_r}} 
\,g_r^{(D)}(x;q,t),$$
the Specialization Formula reads
$$g_r^{(D)}(t^{n-1},t^{n-2},\ldots,1 ;\ q,t)=
t^{r(1-n)} \ \frac{(t^n;q)_r}{(q;q)_r}
\frac{(t^{2n-2};q)_r}{(t^{n-1};q)_r}.$$
We have easily
$g_r^{(D)}(q,t)=g_r^{(C)}(q,t,1)$ and Lemma 3 can be written as 
follows.
\proclaim{Lemma 4} Conjecture 1 yields the specialization
$$g_r^{(D)}(t^{n-1} a,t^{n-2} a,\ldots, a ;\ q,t) = a^r
\frac{(t^n;q)_r}{(q;q)_r}\ {}_2{\phi}_1 \left[\matrix  
q^{-r},\ t^{n-1} \\\ t^{1-n}q^{1-r} \endmatrix;q,t^{2-2n}q/a^2 \right].$$
\endproclaim

As before, for $a=1$ the 
previous expression is a $q$--Chu--Vandermonde sum, and we 
recover the Specialization Formula.

\proclaim{Conjecture 2} For any positive integer $r$ we have
$$g_r^{(D)}(q,t)=\sum_{i=0}^{[r/2]}
g_{r-2i}^{(C)}(q,t,T)\ T^i\ \frac{(1/T;q)_i}{(q;q)_i}\
\frac{(t^nq^{r-2i};q)_{2i}}{(t^{n-1}q^{r-i};q)_i\ (Tt^{n-1}q^{r-2i+1};q)_{i}}.$$
Conversely
$$g_r^{(C)}(q,t,T)=\sum_{i=0}^{[r/2]}
g_{r-2i}^{(D)}(q,t)\ \frac{(T;q)_i}{(q;q)_i}\ 
\frac{(t^nq^{r-2i};q)_{2i}}{(Tt^{n-1}q^{r-i};q)_{i}\ (t^{n-1}q^{r-2i+1};q)_{i}}.$$
Namely, the transition matrix from $g^{(D)}(q,t)$ to 
$g^{(C)}(q,t,T)$ is $\Cal M (t^n,t;1/t,T/t)$ and 
its inverse is $\Cal M (t^n,t;T/t,1/t)$.
\endproclaim

Conjectures $1$ and $2$ are consistent since 
the former, written for $T=t$ and using Theorem $4$, and the latter, 
written for $T=1$, both yield
\proclaim{Conjecture 3} For any positive integer $r$ we have
$$g_r^{(D)}(q,t)=\sum_{i=0}^{[r/2]}
G_{r-2i}(q,t)\ t^i\ \frac{(1/t;q)_i}{(q;q)_i}\
\frac{(t^nq^{r-i};q)_i}{(t^{n-1}q^{r-i};q)_i}\
\frac{1-t^nq^{r-2i}}{1-t^nq^{r-i}}.$$
Conversely
$$G_r(q,t)=\sum_{i=0}^{[r/2]}
g_{r-2i}^{(D)}(q,t)\ \frac{(t;q)_i}{(q;q)_i}\
\frac{(t^nq^{r-2i};q)_i}{(t^{n-1}q^{r-2i+1};q)_i}.$$
Equivalently, the transition matrix from $g^{(D)}(q,t)$ to 
$G(q,t)$ is $\Cal M (t^n,t;1/t,1)$, and 
its inverse is $\Cal M (t^n,t;1,1/t)$.
\endproclaim

The corollary of Lemma 3 shows that Conjecture 3 is true for $x_i=t^{n-i}\ (1\le i\le n)$. Specialization yields another consistency argument.
\proclaim{Lemma 5} Assuming Lemmas 3 and 4, Conjecture 2 is 
true for $x_i=t^{n-i} a$, $ 1\le i\le n$.
\endproclaim
\demo{Proof} A straightforward application of Theorem $1$ with
$u=t^n$, $v=t/a^2$, $x=1/t$, and $y=T/t$.\qed
\enddemo

\head \bf 5. $B$ versus $D$ \endhead

For the type $B$ root system, the set of positive roots is the union 
of $R_{1}=\{\ep_{i}\pm \ep_{j}, 1\le i < j\le n\}$ and 
$R_{2}=\{\ep_{i}, 1\le i\le n\}$. Elements of each 
set have the same length, and we write $t_{\al}=t$ for $\al \in R_{1}$ 
and $t_{\al}=T$ for $\al \in R_{2}$. 

The Weyl group $W$ is the semi-direct product of the 
permutation group $S_{n}$ by $(\Z/2\Z)^n$. 
It acts on $V$ by signed permutation of components. 
The fundamental weights are given by 
$\om_{i}=\sum_{j=1}^i\ep_j$ for $1\le i\le n-1$, and 
$\om_{n}=\frac{1}{2}(\ep_1+\ldots+\ep_n)$. This is the only 
minuscule weight.

The dominant weights $\la\in P^+$ can be identified with vectors 
$\la=\sum_{i=1}^n \la_i \ep_i$, whose components are all integers or all half-integers, and subject to the condition
$\la_1 \ge \la_2\ge \ldots\ge\la_n\ge0$. 
The partial order $\la\geq\mu$ is given by
$\sum_{i=1}^j (\la_i-\mu_i) \in \N$, for $j=1,\ldots, n$.
Writing $t=q^{k}$ and $T=q^{K}$, we have 
$$\rh_k=\suml_{i=1}^n \left((n-i)k+K/2 \right) \ep_i,\quad \quad
  \rh_k^\ast= \suml_{i=1}^n \left((n-i)k+K \right) \ep_i.$$

If we define
$$P_{r\om_1}^{(B)}(x;q,t,T)=\frac{{(q;q)_r}}{{(t;q)_r}} 
\,g_r^{(B)}(x;q,t,T),$$
the Specialization Formula reads
$$g_r^{(B)}(t^{n-1}T,t^{n-2}T,\ldots,T ;\ q,t,T)=
\frac{t^{r(1-n)}}{T^r} \ \frac{(t^n;q)_r}{(q;q)_r}
\ \frac{(Tt^{2n-2};q)_{r}}{(Tt^{2n-2};q)_{2r}}
\ \frac{(T^2t^{2n-2};q)_{2r}}{(Tt^{n-1};q)_r}.$$

\proclaim{Conjecture 4} For any positive integer $r$ we have
$$g_r^{(B)}(q,t,T)=\sum_{i=0}^{r}
g_{r-i}^{(D)}(q,t)\ \frac{(T;q)_i}{(q;q)_i}
\frac{(t^nq^{r-i};q)_i}{(t^{n-1}q^{r-i+1};q)_i}\
\frac{(t^{2n-2}q^{2r-i+1};q)_i}{(Tt^{2n-2}q^{2r-i};q)_i}.$$
Conversely
$$\multline
g_r^{(D)}(q,t)=\sum_{i=0}^{r}
g_{r-i}^{(B)}(q,t,T)\ T^i\ \frac{(1/T;q)_i}{(q;q)_i}\
\frac{(t^nq^{r-i};q)_i}{(t^{n-1}q^{r-i+1};q)_i}\\\times
\frac{(t^{2n-2}q^{2r-2i+1};q)_i}{(Tt^{2n-2}q^{2r-2i+1};q)_i}\
\frac{1-t^{2n-2}q^{2r}}{1-t^{2n-2}q^{2r-i}}.
\endmultline$$
Namely, the transition matrix from $g^{(B)}(q,t,T)$ to  
$g^{(D)}(q,t)$ is $\Cal N (t^n,t^{n-1};T,1)$, and 
its inverse is $\Cal N (t^n,t^{n-1};1,T)$.
\endproclaim
Of course for $T=1$ we recover $g_r^{(D)}(q,t)=g_r^{(B)}(q,t,1)$.

\proclaim{Lemma 6} Conjectures 1 and 4 yield the specialization
$$\multline
g_r^{(B)}(t^{n-1} a,t^{n-2} a,\ldots, a ;\ q,t,T)\\=
(-1)^r t^{r(1-n)}\ \frac{(t^n;q)_r}{(q;q)_r}\
\frac{(Tt^{2n-2};q)_r}{(Tt^{2n-2};q)_{2r}}
\frac{(T^2t^{2n-2};q^2)_{r}(t^{2n-2}q;q^2)_{r}}{(Tt^{n-1};q)_r}
\\\times
{}_4{\phi}_3 \left[\matrix q^{-r},Tt^{2n-2}q^r,-at^{n-1},-1/a\\
t^{n-1}q^{\frac{1}{2}},-t^{n-1}q^{\frac{1}{2}},-Tt^{n-1}\endmatrix;q,q \right].
\endmultline$$
\endproclaim
\demo{Proof} A straightforward application of Lemma 4 and Theorem $2$ with $u=t^n$, $v=t^{n-1}$ and $x=T$.\qed
\enddemo

\proclaim{Corollary} This property is true for $a=T$.
\endproclaim
\demo{Proof} Keeping the notations $u=t^n$, $v=t^{n-1}$ and $x=T$, 
we now have $x=a$. Applying the $q$--Saalsch\"{u}tz formula \cite{3, (II.12)}, we have
$$\align
{}_4{\phi}_3 \left[\matrix q^{-r},av^2q^r,-av,-1/a\\
vq^{\frac{1}{2}},-vq^{\frac{1}{2}},-av\endmatrix;q,q \right]
&={}_3{\phi}_2 \left[\matrix q^{-r},av^2q^r,-1/a\\
vq^{\frac{1}{2}},-vq^{\frac{1}{2}}\endmatrix;q,q \right]\\ 
&= \frac{(-avq^{\frac{1}{2}}q)_r}{(vq^{\frac{1}{2}};q)_r}\
\frac{(q^{\frac{1}{2}}q^{-r}/av;q)_r}
{(-q^{\frac{1}{2}}q^{-r}/v;q)_r}.
\endalign$$
This can be written
$$\align
{}_4{\phi}_3 \left[\matrix q^{-r},av^2q^r,-av,-1/a\\
vq^{\frac{1}{2}},-vq^{\frac{1}{2}},-av\endmatrix;q,q \right]
&= (-1/a)^r \frac{(avq^{\frac{1}{2}};q)_r}
{(vq^{\frac{1}{2}};q)_r}\
\frac{(-avq^{\frac{1}{2}};q)_r}{(-vq^{\frac{1}{2}};q)_r}\\
&= (-1/a)^r \frac{(a^2v^2;q)_{2r}}{(a^2v^2;q^2)_r(v^2q;q^2)_r}.
\endalign$$
By substitution we recover the Specialization Formula. \qed
\enddemo

\head \bf 6. $B$ versus $C$ \endhead

The development of $g_r^{(B)}(q,t,U)$ in 
terms of $g_r^{(C)}(q,t,T)$ (and conversely) can be
immediately written by composing the previous conjectures. 
We did not find a more compact expansion. 

\proclaim{Conjecture 5} For any positive integer $r$ we have
$$\multline
g_r^{(B)}(q,t,U)=\sum_{0\le i+2j\le r}
g_{r-i-2j}^{(C)}(q,t,T)\ T^j\ \frac{(1/T;q)_j}{(q;q)_j}\
\frac{(U;q)_i}{(q;q)_i}\\\times
\frac{(t^nq^{r-i-2j};q)_{i+2j}}{(t^{n-1}q^{r-i-j};q)_{i+j}\ 
(Tt^{n-1}q^{r-i-2j+1};q)_{j}} \
\frac{1-t^{n-1}q^{r-i}}{1-t^{n-1}q^r}\
\frac{(t^{2n-2}q^{2r-i+1};q)_i}{(Ut^{2n-2}q^{2r-i};q)_i}.
\endmultline$$
Conversely
$$\multline
g_r^{(C)}(q,t,T)=\sum_{0\le 2i+j\le r}
g_{r-2i-j}^{(B)}(q,t,U)\ U^j\ \frac{(1/U;q)_j}{(q;q)_j} \
\frac{(T;q)_i}{(q;q)_i} \\
\frac{(t^{2n-2}q^{2r-4i-2j+1};q)_j}{(Ut^{2n-2}q^{2r-4i-2j+1};q)_j}\
\frac{1-t^{2n-2}q^{2r-4i}}{1-t^{2n-2}q^{2r-4i-j}}
 \frac{(t^nq^{r-2i-j};q)_{2i+j}}{(Tt^{n-1}q^{r-i};q)_{i} 
(t^{n-1}q^{r-2i-j+1};q)_{i+j}}.
\endmultline$$
\endproclaim

\head \bf 7. Open problems \endhead

The statements of Lemmas 3, 4 and 6 can be written in a rather similar form by applying the following property.
\proclaim{Lemma 7} We have
$${}_2{\phi}_1 
\left[\matrix q^{-r},\ v\\q^{1-r}/v 
\endmatrix;q,q/v^2a^2 \right]=
(av)^{-r} \frac{(v^2;q)_{r}}{(v;q)_{r}}
{}_4{\phi}_3 \left[\matrix q^{-r},v^2q^{r},av,1/a\\
vq^{\frac{1}{2}},-vq^{\frac{1}{2}},-v\endmatrix;q,q \right].$$
\endproclaim
\demo{Proof} A straightforward consequence of \cite{3, (7.4.12) and (7.4.13)}.\qed
\enddemo

Writing $g_r^{(R)}(a ; q,t,T)$ for
$g_r^{(R)}(t^{n-1} a,t^{n-2} a,\ldots, a ; q,t,T)$, we then have
$$\align
g_r^{(C)}(a ;q,t,T) = (-1)^r &g_r^{(C)}(T^{\frac{1}{2}} ;q,t,T)\\
\times &{}_4{\phi}_3 \left[\matrix q^{-r},T^2t^{2n-2}q^r,-aT^{\frac{1}{2}}t^{n-1},-T^{\frac{1}{2}}/a\\
Tt^{n-1}q^{\frac{1}{2}},-Tt^{n-1}q^{\frac{1}{2}},-Tt^{n-1}\endmatrix;q,q \right]\\
g_r^{(D)}(a ;q,t) = (-1)^r &g_r^{(D)}(1;q,t)\\\times 
&{}_4{\phi}_3 \left[\matrix q^{-r},t^{2n-2}q^r,-at^{n-1},-1/a\\
t^{n-1}q^{\frac{1}{2}},-t^{n-1}q^{\frac{1}{2}},-t^{n-1}\endmatrix;q,q \right]\\
g_r^{(B)}(a ;q,t,T) =
(-T)^r &g_r^{(B)}(T;q,t,T)
\frac{(t^{2n-2}q;q^2)_r}{(T^2t^{2n-2}q;q^2)_r}\\\times
&{}_4{\phi}_3 \left[\matrix q^{-r},Tt^{2n-2}q^r,-at^{n-1},-1/a\\t^{n-1}q^{\frac{1}{2}},-t^{n-1}q^{\frac{1}{2}},-Tt^{n-1}\endmatrix;q,q \right].
\endalign$$

These formulas seem difficult to unify in a general conjecture written in terms of the root system $R$.

We are also in lack of a conjecture for the generating function of $g_r^{(R)}$, except when $R=C$ and $T=t$ (see Theorem 4).

In another paper, we shall present conjectures giving 
for $R \in \{B,C,D\}$,
\smallskip\noindent
(i) a ``generalized Pieri formula'' expanding $g_r^{(R)}g_s^{(R)}$ in terms of the Macdonald polynomials $P_{\la_1\om_1+\la_2\om_2}^{(R)}$,

\smallskip\noindent
(ii) conversely, the expansion of any Macdonald 
polynomial $P_{\la_1\om_1+\la_2\om_2}^{(R)}$ in terms of products 
$g_r^{(R)}g_s^{(R)}$ (``inverse Pieri formula'').

\head \bf 8. $\la$ - rings \endhead

This section and the following will be devoted to the 
proof of Theorem~$4$. This will be done in the language of 
$\la$ - rings, which turns out to be the most efficient. 
Here we only intend to give a short survey of this theory. 
Details and other applications may be found, for instance,
in \cite{6, 7}, and in 
some examples of \cite{10} (see pp. 25, 43, 65 and 79).

The basic idea of the theory of $\la$ - rings is the 
following. A symmetric function $f$ is usually understood as 
\it evaluated \rm on a set of variables $A=\{a_1,a_2,a_3,\ldots\}$, 
the value being denoted $f(A)$. When using $\la$ - rings, this 
interpretation is not the main one. Symmetric 
functions are first understood as \it operators \rm on 
polynomials. Thus any symmetric function $f$ is first 
understood as \it acting on \rm the 
polynomial $P$, mapping $P$ to $f[P]$. Of course the 
standard interpretation may be recovered as a special case.
These statements may be made more precise as follows.

Let $A=\{a_1,a_2,a_3,\ldots\}$ be a (finite or infinite) set of
independent indeterminates, called an alphabet. 
We introduce the generating functions
$$
E_u(A) = \prod_{a\in A} (1+ua), \quad
H_u(A) =  \prod_{a\in A}  \frac{1}{1-ua}, \quad
P_u(A) = \sum_{a\in A} \frac{a}{1-ua},$$
whose development defines symmetric functions known as elementary functions $e_k(A)$, complete functions $h_{k}(A)$, and 
power sums $p_k(A)$, respectively. Each of these three sets generate the 
symmetric algebra $\S(A)$.

We define an action, denoted $[ \ ]$, of $\S(A)$ on the ring 
$\R[A]$ of polynomials in $A$ with real coefficients. Since the power sums $p_{k}$ 
algebraically generate $\S(A)$, it is enough to define the 
action of $p_{k}$ on $\R[A]$. Writing any polynomial as 
$\sum_{c,P} c P$, with $c$ a real constant and $P$ a monomial 
in $(a_{1},a_{2},a_{3},\ldots)$, we define
$$p_{k} \left[\sum_{c,P} c P\right]=\sum_{c,P} c P^{k}.$$
This action extends to $\S[A]$. For instance we obtain
$$ E_u \left[\sum_{c,P} c P\right]= \prod_{c,P} (1+u P)^c \quad , \quad
H_u \left[\sum_{c,P} c P\right]= \prod_{c,P} (1-u P)^{-c}.$$

More generally, we can define an action of $\S(A)$ on the ring of 
rational functions, and even on the ring of formal series, by writing
$$ p_{k} \left( \frac{\sum c P}{\sum d Q}\right) = 
\frac{\sum c P^k}{\sum d Q^k}  $$
with $c,d$ real constants and $P,Q$ monomials in 
$(a_{1},a_{2},a_{3},\ldots)$. This action still extends to $\S(A)$.

As an example, being given a (finite or infinite) alphabet
$A=\{a_1,a_2,a_3,\ldots\}$, let us compute 
$H_u[h_{1}(A)]$ and $H_u[h_{2}(A)]$. We obtain
\cite{10, Example 1.5.10, p. 79}
$$\align
H_u\big[h_1(A)\big]&=\prod_{i}(1-ua_i)^{-1},\\ 
H_u\big[h_2(A)\big]&=\prod_{i} (1-ua_i^2)^{-1}\,
\prod_{i < j} (1-ua_ia_j)^{-1}.
\endalign$$
In the following we shall write 
$$A^{\dag}=h_1(A)=\sum_{i} a_i.$$ 
Bu definition we have
$p_{k} [A^{\dag}]=\sum_{i} a_i^k$. Thus $p_{k} 
[A^{\dag}]=p_{k}(A)$, which yields that for any symmetric function $f$, 
we have $f[A^\dag]=f(A)$.
In particular
$$f(1,q,q^2,q^3,\ldots,q^{m-1})= f 
\left[\frac{1-q^m}{1-q}\right],
\quad f(1,q,q^2,q^3,\ldots)=f\left[\frac{1}{1-q}\right].$$

The following relations are straightforward 
consequences of the previous definitions. For any formal 
series $P,Q$, we have
$$h_r[P+Q]= \sum_{k=0}^r h_{r-k} [P] \ h_k [Q], \quad
e_r[P+Q]= \sum_{k=0}^r e_{r-k} [P] \ e_k [Q].$$
Or equivalently
$$\align
&H_u[P+Q]=H_u[P]\,H_u[Q],\quad \quad E_u[P+Q]=E_u[P]\,E_u[Q]\\
&H_u[P-Q]=H_u[P]\,{H_u[Q]}^{-1},\quad E_u[P-Q]=E_u[P]\,{E_u[Q]}^{-1}.
\endalign$$
As an application, for a finite alphabet 
$A=\{a_1,a_2,\ldots,a_m\}$ we may write
$$\align
H_1 \left[\frac{uA^\dag}{1-q}\right]&=\prod_{i\geq0} 
H_1 \big[uq^iA^\dag\big] 
=\prod_{k=1}^{m} \ \prod_{i\geq0} H_1 \big[uq^ia_k\big]\\
&=\prod_{k=1}^{m} \ \prod_{i\geq0} \frac{1}{1-uq^ia_k}
=\prod_{k=1}^{m} \frac{1}{{(ua_k;q)}_{\infty}},
\endalign$$
and
$$H_1 \left[\frac{1-t}{1-q}A^\dag\right]= 
H_1 \left[\frac{A^\dag}{1-q}\right] 
{\left(H_1 \left[\frac{tA^\dag}{1-q}\right]\right)}^{-1}.$$
Finally we obtain
$$H_1\left[\frac{1-t}{1-q}\,A^\dag\right]= \sum_{r\ge 0} 
h_{r}\left[\frac{1-t}{1-q}\,A^\dag\right]
= \prod_{i=1}^{m} \frac{{(ta_i;q)}_{\infty}}{{(a_i;q)}_{\infty}}.$$

\head \bf 9. Proof of Theorem 4 \endhead

In this section we assume that $R$ is of type $C_n$ with 
$T=t$. We define $x_i=e^{\ep_i}$ for $i=1\ldots n$, and regard 
elements of $\A$ as Laurent polynomials of $n$ variables $x_1,\ldots,x_n$.

The dual root system $R^\lor=B_n$ has one minuscule 
weight $\pi=\frac{1}{2}(\ep_1+\ldots+\ep_n)$. 
With the notations of Section $1$, we have
$$\Phi_\pi =
\prod_{i=1}^n\frac {1-tx_i^2} {1-x_i^2} \
\prod_{1\le i < j \le n}\frac {1-tx_ix_j} {1-x_ix_j},$$ 
and the translation operator $T_\pi$ acts on $\A$ by
$$T_{\pi }f(x_1,\ldots,x_n)=f(q^{\frac{1}{2}}x_1,\ldots,q^{\frac{1}{2}}x_n).$$

The Weyl group $W$ is the semi-direct product of the 
permutation group $S_{n}$ by $(\Z/2\Z)^n$. 
It acts on $V$ by signed permutation of components. 
Hence the $W$-orbit of $\pi$ is formed by vectors
$\frac{1}{2}(\ta_1\ep_1+\ldots+\ta_n\ep_n)$ with 
$\ta\in(-1,+1)^n$.

The Macdonald operator $E_\pi$ can be written as
$$E_\pi f=\sum_{\ta\in(-1,+1)^n} 
\prod_{i=1}^n\frac {1-tx_{i}^{2\ta_i}} 
{1-x_{i}^{2\ta_i}}\
\prod_{1\le i < j \le n}
\frac {1-t x_{i}^{\ta_i}x_{j}^{\ta_j}} 
{1-x_{i}^{\ta_i}x_{j}^{\ta_j}}\ 
f(q^{\ta_1/2}x_1,\ldots,q^{\ta_n/2}x_n).$$
Up to a constant, the Macdonald polynomial 
$P_{r\om_1}$ is defined by
$$E_\pi P_{r\om_1}=\prod_{i=1}^{n-1}(t^i+1) \
(t^nq^{r/2}+q^{-r/2}) \,P_{r\om_1}.$$
This normalization constant is given by the condition
$$P_{r\om_1} = m_{r\om_1}+ \text{lower terms},$$
where the orbit-sum $m_{r\om_1}$ is given by 
$m_{r\om_1}(x)=\sum_{i=1}^n (x_i^r+1/x_i^r)$.

The following notations will be used till the end of this 
paper. We write 
$$X_{+}=\{x_1,\ldots,x_n\},\quad 
X_-=\{1/x_1,\ldots,1/x_n\}, \quad X=X_{+} \cup X_-,$$
so that $X^{\dag}=\sum_{i=1}^n (x_i+1/x_i)$.
Observe that
$\Phi_\pi =H_1\big[(1-t)h_2(X_{+})\big]$.

We consider the generating series
$$H_1\left[\frac{1-t}{1-q}X^{\dag}\right]=
\sum_{r\ge 0} u^r\ h_r\left[\frac{1-t}{1-q}X^{\dag}\right]=
\prod_{i=1}^n \frac{(tux_i;q)_\infty} {(ux_i;q)_\infty}\
\frac{(tu/x_i;q)_\infty} {(u/x_i;q)_\infty}.$$
Theorem $4$ states that
$$g_{r}^{(C)}(x;q,t,t)=h_r\left[\frac{1-t}{1-q}X^{\dag}\right].$$
Since it is well known (\cite{10, p. 314}, \cite{7, p. 237}) that
for any alphabet $A=\{a_1,\ldots,a_m\}$ one has
$$h_r\left[\frac{1-t}{1-q}A^{\dag}\right]=
\frac{(t;q)_r}{(q;q)_r} \ \sum_{i=1}^m a_i^r + \text{other terms},$$
we have only to prove
$$E_\pi \ h_r\left[\frac{1-t}{1-q}X^{\dag}\right]
=\prod_{i=1}^{n-1}(t^i+1) \
(t^nq^{r/2}+q^{-r/2}) \ h_r\left[\frac{1-t}{1-q}X^{\dag}\right].$$
This will be established under the following equivalent form.
\proclaim{Theorem 4'}
We have 
$$E_\pi H_1\left[\frac{1-t}{1-q}X^{\dag}\right]=
\prod_{i=1}^{n-1}(t^i+1) \left(t^n 
H_1\left[\frac{1-t}{1-q}q^{\frac{1}{2}}X^{\dag}\right]
+H_1\left[\frac{1-t}{1-q}q^{-\frac{1}{2}}X^{\dag}\right]\right).$$
\endproclaim

In order to evaluate the left-hand side of this 
identity, we shall need the following trick. For $\ta=\pm 1$ we have
$$xq^{\ta/2}+\frac{1}{xq^{\ta/2}}=
q^{\frac{1}{2}}\left(x+\frac{1}{x}\right)+q^{-\frac{1}{2}}(1-q)x^{-\ta},$$
which is checked separately for $\ta=1$ and $\ta=-1$.
Consequently we write
$$\align
H_1\left[\frac{1-t}{1-q}
\sum_{i=1}^n 
\left(x_iq^{\ta_i/2}+\frac{1}{x_iq^{\ta_i/2}}\right)\right]&=
H_1\left[\frac{1-t}{1-q}q^{\frac{1}{2}}X^{\dag}+
(1-t)q^{-\frac{1}{2}}\sum_{i=1}^n x_i^{-\ta_i}\right]\\
&=H_1\left[(1-t)q^{-\frac{1}{2}}\sum_{i=1}^n x_i^{-\ta_i}\right] 
\, H_1\left[\frac{1-t}{1-q}q^{\frac{1}{2}}X^{\dag}\right]\\
&=\prod_{i=1}^n \frac{1-tq^{-\frac{1}{2}}x_i^{-\ta_i}}
{1-q^{-\frac{1}{2}}x_i^{-\ta_i}}\
H_1\left[\frac{1-t}{1-q}q^{\frac{1}{2}}X^{\dag}\right].
\endalign$$
Similarly, on the right-hand side we get
$$\align
H_1\left[\frac{1-t}{1-q}q^{-\frac{1}{2}}X^{\dag}\right]&=
H_1\left[\frac{1-t}{1-q}q^{\frac{1}{2}}X^{\dag}+(1-t)q^{-\frac{1}{2}}X^{\dag}\right]\\
&=H_1\big[(1-t)q^{-\frac{1}{2}}X^{\dag}\big]\
H_1\left[\frac{1-t}{1-q}q^{\frac{1}{2}}X^{\dag}\right] \\
&=\prod_{i=1}^n\frac{1-tq^{-\frac{1}{2}}x_i}{1-q^{-\frac{1}{2}}x_i}\,
\frac{1-tq^{-\frac{1}{2}}/x_i}{1-q^{-\frac{1}{2}}/x_i}\
H_1\left[\frac{1-t}{1-q}q^{\frac{1}{2}}X^{\dag}\right].
\endalign$$
By substitution, we see that the statement of Theorem $4'$ is 
equivalent to the following rational identity, 
written with $u=q^{-\frac{1}{2}}$,
$$\multline
\sum_{\ta\in(-1,+1)^n} 
\prod_{i=1}^n\frac {1-tx_{i}^{2\ta_i}} 
{1-x_{i}^{2\ta_i}}\
\frac{1-tux_i^{-\ta_i}}
{1-ux_i^{-\ta_i}}\
\prod_{1\le i < j \le n}
\frac {1-tx_{i}^{\ta_i}x_{j}^{\ta_j}} 
{1-x_{i}^{\ta_i}x_{j}^{\ta_j}}=\\
\prod_{i=1}^{n-1}(t^i+1) \left(t^n+
\prod_{i=1}^n \frac{1-tux_i}
{1-ux_i}\, \frac{1-tu/x_i}
{1-u/x_i}\right).
\endmultline$$

Writing $T_{i}$ for the operator $x_i\rightarrow1/x_i$, 
we are led to prove Theorem $4$ under the following equivalent form.
\proclaim{Theorem 5}
We have 
$$\multline
(1+T_1)\cdots(1+T_n)\left(
\prod_{i=1}^n\frac {1-tx_i^2} {1-x_i^2}\,
\frac{1-tu/x_i}{1-u/x_i} \
\prod_{1\le i < j \le n}
\frac {1-tx_ix_j} {1-x_ix_j}\right)=\\
\prod_{i=1}^{n-1}(t^i+1) \left(t^n+
\prod_{i=1}^n \frac{1-tux_i}{1-ux_i}\ 
\frac{1-tu/x_i}{1-u/x_i}\right).
\endmultline$$
\endproclaim

Both sides of this identity are rational 
functions of $u$ having poles at $u=x_i$ and $u=1/x_i$, for 
$i=1,\ldots,n$. We first prove 
that their constant terms are equal, i.e. that the statement 
is true for $u=0$.
\proclaim{Lemma 8} We have
$$(1+T_1)\cdots(1+T_n)\ \Phi_\pi=\prod_{i=1}^{n}(t^i+1)$$
\endproclaim
\demo{Proof} This is a direct consequence of Weyl's denominator formula
$$\delta = e^{-\rho}\prod_{\alpha\in R^+} (e^\alpha -1)= \sum_{w\in W} \det
(w)\, e^{w\rho},$$
with $\rh=\sum_{i=1}^n (n-i+1)\ep_i$, from which follows
$$\Phi_\pi= \delta^{-1}\sum_{w\in W} \det (w)\, 
t^{\langle \pi, \rho+w\rho\rangle} e^{w\rho}.$$
Since for any $w\in W$, we have $w\delta = \det(w) \delta$, we obtain
$$(1+T_1)\cdots(1+T_n)\ \Phi_\pi=  
\sum_{\tau\in W(\pi)} t^{\langle \pi+\tau, \rho \rangle}
=\sum_{\ta\in(-1,+1)^n} \prod_{i=1}^n t^{\frac{1}{2}(n-i+1)(1+\ta_i)}.\qed$$
\enddemo

\demo{Proof of Theorem 5} It is sufficient to prove that both sides of the 
identity have the same residue at each of their poles, i.e. 
at $u=x_i$ and $u=1/x_i$, $1\le i\le n$. 
By symmetry, this has only to be checked for some 
$x_i$, say $x_n$. We shall only do it at $u=x_n$, the proof at 
$u=1/x_n$ being similar.

If $A=\{a_1,\ldots,a_m\}$ is an arbitrary alphabet, 
we have 
$$\prod_{i=1}^{m} \frac{tu-a_i}{u-a_i}=t^m+
(t-1)\sum_{i=1}^{m}\frac{a_i}{u-a_i} \
\underset{j\neq i}\to{\prod_{j=1}^m}
\frac{ta_i-a_j}{a_i-a_j}.$$
This decomposition as a sum of partial fractions is actually 
a Lagrange interpolation (see \cite{5}, and also \cite{7, p. 236}). 
We first apply it to the right-hand side of the identity. Its residue 
at $u=x_n$ is given by
$$\prod_{i=1}^{n-1}(t^i+1)\ x_n\ (t-1)\ \prod_{j=1}^{n-1}
\frac{tx_n-x_j}{x_n-x_j} 
\prod_{j=1}^{n}\frac{tx_n-1/x_j}{x_n-1/x_j}.$$
We then apply the Lagrange interpolation to
$$\prod_{i=1}^n \frac{tu-x_i^{\ta_i}}{u-x_i^{\ta_i}},$$
on the left-hand side of the identity.
Only fractions with $\ta_n=1$ contribute to the residue 
at $u=x_n$, which can be written as
$$\multline
x_{n}(t-1)\frac {1-tx_{n}^{2}} 
{1-x_{n}^{2}}\\\times \sum_{\ta\in(-1,+1)^{n-1}} 
\prod_{i=1}^{n-1}\frac {1-tx_{i}^{2\ta_i}} 
{1-x_{i}^{2\ta_i}}\
\frac{1-tx_{n}x_i^{-\ta_i}}
{1-x_{n}x_i^{-\ta_i}}\
\frac {1-tx_{i}^{\ta_i}x_{n}} 
{1-x_{i}^{\ta_i}x_{n}}\
\prod_{1\le i < j \le n-1}
\frac {1-tx_{i}^{\ta_i}x_{j}^{\ta_j}} 
{1-x_{i}^{\ta_i}x_{j}^{\ta_j}}.
\endmultline$$

By identification of residues on both sides, and using an obvious induction on 
$n$, we are led to prove the identity
$$\multline
(1+T_1)\cdots(1+T_n)\left(
\prod_{i=1}^n\frac {1-tx_i^2} {1-x_i^2}\
\frac{1-tux_i}{1-ux_i}\ \frac{1-tu/x_i}{1-u/x_i}
\prod_{1\le i < j \le n}
\frac {1-tx_ix_j} {1-x_ix_j}\right)=\\
\prod_{i=1}^{n}(t^i+1) \ \frac{1-tux_i}{1-ux_i}\ 
\frac{1-tu/x_i}{1-u/x_i}.
\endmultline$$
Since
$$\prod_{i=1}^n \frac{1-tux_i}{1-ux_i}\ 
\frac{1-tu/x_i}{1-u/x_i}=H_1\big[(1-t)uX^{\dag}\big]$$
is obviously 
invariant under any $T_i$, the left-hand side may be written
$$\prod_{i=1}^n \frac{1-tux_i}{1-ux_i}\ 
\frac{1-tu/x_i}{1-u/x_i}\ (1+T_1)\cdots(1+T_n) \Phi_\pi.$$
We conclude by applying Lemma 8.\qed
\enddemo

\head \bf 10. Proofs of Theorems 1, 2 and 3 \endhead

We shall implicitly use many of the formulas about $q$-shifted 
factorials, listed in Appendix I of~\cite{3}. In particular we shall write
$$(a_1,a_2,\dots,a_{r};q)_{i}=(a_1;q)_i\ (a_2;q)_i\ \dots \ 
(a_r;q)_i.$$

\proclaim{Theorem 1}We have the following transformation 
between ${}_2{\phi}_1$ series
$$\multline
\frac{(u;q)_r}{(q;q)_r}\ {}_2{\phi}_1 
\left[\matrix q^{-r}, \ ux\\q^{1-r}/ux \endmatrix;q,qv/u^2x \right]\\=
\sum_{i=0}^{[r/2]} 
\frac{(u;q)_{r-2i}}{(q;q)_{r-2i}}\ 
{}_2{\phi}_1 
\left[\matrix q^{2i-r}, \ uy\\q^{1+2i-r}/uy \endmatrix;q,qv/u^2y \right]\ y^iv^i\ 
\frac{(x/y;q)_i\ (uq^{r-2i};q)_{2i}}{(q,uxq^{r-i},uyq^{r-2i+1};q)_{i}}.
\endmultline$$
\endproclaim
\demo{Proof~\cite{14}} With $a=v/u$, applying \cite{3, (3.4.7)}, we have
$$\multline\frac{(u;q)_r}{(q;q)_r}\ 
{}_2{\phi}_1 \left[\matrix  q^{-r},\ ux \\
\ q^{1-r}/ux \endmatrix;q,qv/u^2x \right]=
a^{r} \ \frac{(u,1/a^2;q)_r}{(q,1/a;q)_r}\\\times
{}_8{\phi}_7 \left[\matrix 
q^{-r}a,q(q^{-r}a)^{\frac{1}{2}},-q(q^{-r}a)^{\frac{1}{2}},aux,
q^{-\frac{r}{2}},-q^{-\frac{r}{2}},
q^{\frac{1}{2}-\frac{r}{2}},-q^{\frac{1}{2}-\frac{r}{2}}\\
(q^{-r}a)^{\frac{1}{2}},-(q^{-r}a)^{\frac{1}{2}},q^{1-r}/ux,
q^{1-\frac{r}{2}}a,-q^{1-\frac{r}{2}}a, 
q^{\frac{1}{2}-\frac{r}{2}}a,-q^{\frac{1}{2}-\frac{r}{2}}a 
\endmatrix;q,qa/ux \right],
\endmultline$$
or equivalently
$$\multline \frac{(u;q)_r}{(q;q)_r}\ 
{}_2{\phi}_1 \left[\matrix  q^{-r},\ ux \\
\ q^{1-r}/ux \endmatrix;q,qv/u^2x \right]=
a^{r}\ \frac{(u,1/a^2;q)_r}{(q,1/a;q)_r}\\\times
\sum_{k=0}^{[r/2]}\frac{1-aq^{2k-r}}{1-aq^{-r}}
\frac{(q^{-r}a,aux;q)_k}{(q^{1-r}/ux,q;q)_k}\ 
\frac{(q^{-r};q)_{2k}}{(q^{1-r}a^2;q)_{2k}}\ (qa/ux)^k.
\endmultline$$
This yields
$$\multline\frac{(u;q)_{r-2i}}{(q;q)_{r-2i}}\ 
{}_2{\phi}_1 \left[\matrix  q^{2i-r},\ uy \\
\ q^{1+2i-r}/uy \endmatrix;q,qv/u^2y \right]\\=
(v/u)^{r}\ \frac{(u,u^2/v^2;q)_r}{(q,u/v;q)_r}\ (q/u)^{2i}\ 
\frac{(q^{-r},q^{1-r}v/u;q)_{2i}}{(q^{1-r}/u,q^{1-r}v^2/u^2;q)_{2i}}\\
\sum_{k=0}^{[r/2]-i}\frac{1-q^{2i+2k-r}v/u}{1-q^{2i-r}v/u}
\frac{(q^{2i-r}v/u,vy;q)_k}{(q^{1+2i-r}/uy,q;q)_k}\ 
\frac{(q^{2i-r};q)_{2k}}{(q^{1+2i-r}v^2/u^2;q)_{2k}}\ (qv/u^2y)^k.
\endmultline$$
If we substitute this value in the right-hand side of the 
identity, and put $k=j-i$, we obtain
$$\multline 
(v/u)^{r}\ \frac{(u,u^2/v^2;q)_r}{(q,u/v;q)_r}\
\sum_{j=0}^{[r/2]} \sum_{i=0}^{j} \
(q/u)^{2i}\ \frac{(q^{-r},q^{1-r}v/u;q)_{2i}}
{(q^{1-r}/u,q^{1-r}v^2/u^2q)_{2i}}\\\times
\frac{1-q^{2j-r}v/u}{1-q^{2i-r}v/u} 
\frac{(q^{2i-r}v/u,vy;q)_{j-i}}{(q^{1+2i-r}/uy,q;q)_{j-i}}\ 
\frac{(q^{2i-r};q)_{2j-2i}}{(q^{1+2i-r}v^2/u^2;q)_{2j-2i}}\ 
(qv/u^2y)^{j-i}\\\times
y^iv^i\ \frac{(x/y;q)_i\ (uq^{r-2i};q)_{2i}}
{(q,uxq^{r-i},uyq^{r-2i+1};q)_{i}}.
\endmultline$$
This may be written as
$$\multline 
(v/u)^{r}\ 
\frac{(u,u^2/v^2;q)_r}{(q,u/v;q)_r}\\\times
\sum_{j=0}^{[r/2]} \frac{1-q^{2j-r}v/u}{1-q^{-r}v/u}\
\frac{(q^{-r}v/u,vy;q)_j}{(q^{1-r}/uy,q;q)_j}\
\frac{(q^{-r};q)_{2j}}{(q^{1-r}v^2/u^2;q)_{2j}}\ (qv/u^2y)^{j} 
\\\times \sum_{i=0}^{j}
\frac{(q^{1-r}/uy;q)_{2i}}{(q^{1-r}/u;q)_{2i}}\ 
\frac{(q^{-j},q^{j-r}v/u,x/y;q)_i\ (uq^{r-2i};q)_{2i}}
{(q^{1-j}/vy,q^{1+j-r}/uy,q,uxq^{r-i},uyq^{r-2i+1};q)_i}\
\ (q^2/vy)^i.
\endmultline$$
The sum over $i$ reads
$$\multline 
\sum_{i=0}^{j}
\frac{1-uyq^{r-2i}}{1-uyq^{r-i}}\ 
\frac{(q^{-j},q^{j-r}v/u,q^{1-r}/uy,x/y;q)_i}
{(q^{1-j}/vy,q^{1+j-r}/uy,q^{1-r}/ux,q;q)_i}\
\ (q^2/vx)^i\\
={}_6{\phi}_5 \left[\matrix 
q^{-r}/uy,q(q^{-r}/uy)^{\frac{1}{2}},-q(q^{-r}/uy)^{\frac{1}{2}},x/y,
q^{j-r}v/u,q^{-j}\\
(q^{-r}/uy)^{\frac{1}{2}},-(q^{-r}/uy)^{\frac{1}{2}},q^{1-r}/ux,
q^{1-j}/vy,q^{j+1-r}/uy, 
\endmatrix;q,q/vx \right].
\endmultline$$
By \cite{3, (II.21)} this terminating very-well-poised $_6\phi_5$ 
sum equals
$$
\frac{(q^{1-r}/uy,q^{1-j}/vx;q)_{j}}{(q^{1-r}/ux,q^{1-j}/vy;q)_j}=
\frac{(q^{1-r}/uy,vx;q)_j}{(q^{1-r}/ux,vy;q)_j}\ (y/x)^j.$$
Finally we have proved that the right-hand side of the identity is
$$\multline 
(v/u)^{r}\ \frac{(u,u^2/v^2;q)_r}{(q,u/v;q)_r}\\\times
\sum_{j=0}^{[r/2]}\frac{1-q^{2j-r}v/u}{1-q^{-r}v/u}
\frac{(q^{-r}v/u,vx;q)_j}{(q^{1-r}/ux,q;q)_j}\ 
\frac{(q^{-r};q)_{2j}}{(q^{1-r}v^2/u^2;q)_{2j}}\ (qv/u^2x)^j.
\endmultline$$
Hence the statement.\qed
\enddemo

\proclaim{Lemma 9} We have
$$\multline
{}_4{\phi}_3 \left[\matrix q^{-n},a^2,qa,b\\
a,qa^2/b,q^{1+n}a^2\endmatrix;q,q^{1+n}a/b 
\right]\\=
\frac{(qa^2,-1;q)_{n}}
{(q^{\frac{1}{2}}a,-q^{\frac{1}{2}}a;q)_{n}}\ 
{}_4{\phi}_3 \left[\matrix 
q^{-n},-qa/b,q^{\frac{1}{2}}a,-q^{\frac{1}{2}}a\\
-q^{1-n},qa^2/b,-qa\endmatrix;q,q \right].$$
\endmultline$$
\endproclaim
\demo{Proof~\cite{14}}The left-hand side obviously equals
$${}_8{\phi}_7 \left[\matrix 
a^2,qa,-qa,b,-a,q^{\frac{1}{2}}a,-q^{\frac{1}{2}}a,q^{-n}\\
a,-a,qa^2/b,-qa,q^{\frac{1}{2}}a,
-q^{\frac{1}{2}}a,q^{1+n}a^2\endmatrix;q,q^{1+n}a/b 
\right].$$
By applying Watson's transformation formula \cite{3, (III.19)}, it can be transformed into the right-hand side. 
\qed
\enddemo

\proclaim{Theorem 2} Defining
$$\frac{(u;q)_{r}}{(q;q)_{r}}\ H_r=
\sum_{i=0}^r a^{r-i} \frac{(u;q)_{r-i}}{(q;q)_{r-i}}\
{}_2{\phi}_1 \left[\matrix q^{i-r},\ v\\q^{1+i-r}/v 
\endmatrix;q,q/v^2a^2 \right]
\frac{(x,uq^{r-i},v^2q^{2r-i+1};q)_i}
{(q,vq^{r-i+1},xv^2q^{2r-i};q)_i},$$
we have
$$H_r=(-1/v)^r\ 
\frac{(vq^{\frac{1}{2}},-vq^{\frac{1}{2}},-xv,xv^2;q)_r}
{(xv^2;q)_{2r}}\
{}_4{\phi}_3 \left[\matrix q^{-r},xv^2q^r,-av,-1/a\\
vq^{\frac{1}{2}},-vq^{\frac{1}{2}},-xv\endmatrix;q,q \right].$$
\endproclaim
\demo{Proof~\cite{14}} Applying Lemma 7 we obtain
$$\multline  H_r= \frac{(q;q)_{r}}{(v;q)_{r}}\
\sum_{i=0}^r \frac{1-vq^{r-i}}{1-vq^{r}}\ 
 \frac{(v^2;q)_{r-i}}{(q;q)_{r-i}}\
\frac{(x,v^2q^{2r-i+1};q)_{i}}{(q,xv^2q^{2r-i};q)_{i}}\ v^{i-r}\\\times
{}_4{\phi}_3 \left[\matrix q^{i-r},v^2q^{r-i},av,1/a\\
vq^{\frac{1}{2}},-vq^{\frac{1}{2}},-v\endmatrix;q,q \right].
\endmultline$$
Reversing the order of summation, this may be rewritten as
$$\multline H_r=v^{-r}\ 
\frac{(v^2;q)_{r}}{(v;q)_{r}}\
\sum_{i=0}^r \frac{(q^{-r},v^2q^{r},av,1/a;q)_{i}}
{(q,vq^{\frac{1}{2}},-vq^{\frac{1}{2}},-v;q)_{i}}\ q^i\\\times
{}_4{\phi}_3 \left[\matrix q^{-2r}/v^2,q^{1-r}/v,x,q^{i-r}\\
q^{-r}/v,q^{1-2r}/xv^2,q^{1-r-i}/v^2\endmatrix;q,q^{1-i}/xv \right].
\endmultline$$
Using Lemma 9 with $a:=q^{-r}/v$, $b:=x$ and $n:=r-i$, the ${}_4{\phi}_3$ 
series equals
$$\multline \frac{(q^{1-2r}/v^2,-1;q)_{r}}
{(q^{\frac{1}{2}-r}/v,-q^{\frac{1}{2}-r}/v;q)_{r}} \
\ \frac{(vq^{\frac{1}{2}},-vq^{\frac{1}{2}};q)_{i}}
{(v^2q^{r},-q^{1-r};q)_{i}} \\\times {}_4{\phi}_3 \left[\matrix 
q^{i-r},q^{\frac{1}{2}-r}/v,-q^{\frac{1}{2}-r}/v,-q^{1-r}/xv\\
q^{1-2r}/xv^2,-q^{1-r}/v,-q^{1-r+i}\endmatrix;q,q \right],
\endmultline$$
which yields
$$\multline H_r= \frac{(v^2,v^2q^r,-1;q)_{r}}
{(v,vq^{\frac{1}{2}},-vq^{\frac{1}{2}};q)_{r}}\ v^{-r} 
\ q^{-{r \choose 2}}\
\sum_{i=0}^r q^i \ 
\frac{(q^{-r},av,1/a;q)_{i}}{(q,-v,-q^{1-r};q)_{i}}
\\\times  {}_4{\phi}_3 \left[\matrix 
q^{i-r},q^{\frac{1}{2}-r}/v,-q^{\frac{1}{2}-r}/v,-q^{1-r}/xv\\
q^{1-2r}/xv^2,-q^{1-r}/v,-q^{1-r+i}\endmatrix;q,q \right].
\endmultline$$
Reversing again the order of summation, this may be rewritten as 
$$\multline H_r= (-v,-1;q)_{r}\ v^{-r} 
\ q^{-{r \choose 2}}\
\sum_{i=0}^r q^i \
\frac{(q^{-r},q^{\frac{1}{2}-r}/v,-q^{\frac{1}{2}-r}/v,-q^{1-r}/xv;q)_{i}}
{(q,q^{1-2r}/xv^2,-q^{1-r}/v,-q^{1-r};q)_{i}}\\\times
{}_3{\phi}_2 \left[\matrix 
q^{i-r},av,1/a\\-v,-q^{1-r+i}\endmatrix;q,q\right].
\endmultline$$
Since it is a $q$-Saalsch\"{u}tz sum \cite{3, (II.12)}, the $_3{\phi}_2$ series equals
$$ \frac{(-av,-1/a;q)_{r-i}}{(-v,-1;q)_{r-i}}=\frac{(-av,-1/a;q)_r}
{(-v,-1;q)_r} \ \frac{(-q^{1-r},-q^{1-r}/v;q)_{i}}
{(-q^{1-r}/av,-aq^{1-r};q)_{i}}.$$
Finally we obtain
$$H_r=(-av,-1/a;q)_r\ v^{-r}\ q^{-{r \choose 2}}\
{}_4{\phi}_3 \left[\matrix 
q^{-r},q^{\frac{1}{2}-r}/v,-q^{\frac{1}{2}-r}/v,-q^{1-r}/xv\\
q^{1-2r}/xv^2,-q^{1-r}/av,-aq^{1-r}\endmatrix;q,q \right],$$
and we conclude easily.\qed
\enddemo

\proclaim{Theorem 3} 
We have
$$\multline
a^r \frac{(u;q)_{r}}{(q;q)_{r}}\
{}_2{\phi}_1 
\left[\matrix q^{-r},\ v\\q^{1-r}/v 
\endmatrix;q,q/v^2a^2 \right]\\
= \sum_{i=0}^r \frac{(u;q)_{r-i}}{(q;q)_{r-i}}\ H_{r-i}
\ x^i\ \frac{(1/x,uq^{r-i},v^2q^{2r-2i+1};q)_{i}}
{(q,vq^{r-i+1},xv^2q^{2r-2i+1};q)_{i}}\
\frac{1-v^2q^{2r}}{1-v^2q^{2r-i}}.
\endmultline$$
\endproclaim
\demo{Proof~\cite{14}} First observe that the right-hand side is
$$\multline
\sum_{i=0}^r \frac{(u;q)_i}{(q;q)_i}\ 
\frac{(1/x, uq^i,v^2q^{2i+1};q)_{r-i}}
{(q,vq^{i+1},xv^2q^{2i+1};q)_{r-i}}\
\frac{1-v^2q^{2r}}{1-v^2q^{r+i}} \ x^{r-i}H_i\\ = x^r
\frac{(u,1/x,v^2;q)_r}{(q,vq,xv^2q;q)_r} \frac{1-v^2q^{2r}}{1-v^2}
\sum_{i=0}^r \frac{(q^{-r},vq,v^2q^r;q)_{i}}
{(q,xq^{1-r},xv^2q^{1+r};q)_{i}}\
\frac{(xv^2q;q)_{2i}}{(v^2q;q)_{2i}}
(-q/v)^i\\ \qquad \qquad \qquad \times
\frac{(vq^{\frac{1}{2}},-vq^{\frac{1}{2}},-xv,xv^2;q)_i}
{(xv^2;q)_{2i}}\
{}_4{\phi}_3 \left[\matrix q^{-i},xv^2q^i,-av,-1/a\\
vq^{\frac{1}{2}},-vq^{\frac{1}{2}},-xv\endmatrix;q,q \right]\\
=x^r \frac{(u,1/x,v^2,-vq;q)_r}{(q,v,-v,xv^2q;q)_r}
\sum_{i=0}^r \frac{1-xv^2q^{2i}}{1-xv^2}
\frac{(q^{-r},v^2q^r,xv^2,-xv;q)_i}
{(q,xq^{1-r},xv^2q^{1+r},-qv;q)_i}
(-q/v)^i\\\times
{}_4{\phi}_3 \left[\matrix q^{-i},xv^2q^i,-av,-1/a\\
vq^{\frac{1}{2}},-vq^{\frac{1}{2}},-xv\endmatrix;q,q \right].
\endmultline$$
Expanding the ${}_4{\phi}_3$ series over the index $j\le i$, putting $i=j+k$, and summing over $k$, we obtain
$$\multline
x^r \frac{(u,1/x,v^2,-vq;q)_r}{(q,v,-v,xv^2q;q)_r}
\sum_{j=0}^r \frac{(xv^2q;q)_{2j}}{(v^2q;q)_{2j}}
\frac{(vq,-av,-1/a,v^2q^r,q^{-r};q)_j}
{(q,xq^{1-r},xv^2q^{1+r};q)_j}
(q/v)^j \ q^{-{j \choose 2}}\\\times
{}_6{\phi}_5 \left[\matrix xv^2q^{2j},q(xv^2q^{2j})^{\frac{1}{2}},-q(xv^2q^{2j})^{\frac{1}{2}},-xvq^j,v^2q^{r+j},q^{j-r}\\ 
(xv^2q^{2j})^{\frac{1}{2}},-(xv^2q^{2j})^{\frac{1}{2}},
-vq^{j+1},xq^{j-r+1},xv^2q^{j+r+1}\endmatrix;q,-q^{1-j}/v \right].
\endmultline$$
By \cite{3, (II.21)} this terminating very-well-poised $_6\phi_5$ 
sum equals
$$\frac{(xv^2q^{2j+1},-q^{1-r}/v;q)_{r-j}}
{(-vq^{j+1},xq^{j-r+1};q)_{r-j}}= \Big(\frac{-1}{xv}\Big)^{r} 
v^j q^{{j \choose 2}} \frac{(xv^2q,-v;q)_r}{(1/x,-vq;q)_r}
\frac{(-vq,xv^2q^{1+r},xq^{1-r};q)_j}{(-v;q)_{j}(xv^2q;q)_{2j}}.$$
Thus, on simplification, the right-hand side may be written 
$$\frac{(u,v^2;q)_r}{(q,v;q)_r}\ (-v)^{-r}
\ {}_4{\phi}_3 \left[\matrix q^{-r},v^2q^r,-av,-1/a\\
vq^{\frac{1}{2}},-vq^{\frac{1}{2}},-v\endmatrix;q,q \right].$$
Applying Lemma 7 with $-a$ substituted to $a$, we obtain the left-hand side.\qed
\enddemo

\head \bf Acknowledgements \endhead

It is a pleasure to thank Alain Lascoux, Mizan Rahman and 
Michael Schlosser for their interest and generous help.

\enddocument